\documentclass[11pt]{article}
\usepackage{amsmath}
\usepackage{amssymb}
\usepackage{amsthm}
\usepackage[usenames]{color}
\usepackage{amscd}
\usepackage{dsfont}

\def\1{{\bf 1}}

\def\lcm{\operatorname{lcm}}

\newtheorem{theorem}{Theorem}
\newtheorem{cor}{Corollary}

\begin{document}

\title{\bf On the number of cyclic subgroups of \\ a finite Abelian group}
\author{L\'aszl\'o T\'oth \thanks{The author gratefully acknowledges support from the Austrian Science
Fund (FWF) under the project Nr. M1376-N18.}
\\ Institute of Mathematics, Department of Integrative
Biology \\ Universit\"at f\"ur Bodenkultur, Gregor Mendel-Stra{\ss}e
33, A-1180 Wien, Austria \\ and \\ Department of Mathematics,
University of P\'ecs \\ Ifj\'us\'ag u. 6, H-7624 P\'ecs, Hungary \\
E-mail: ltoth@gamma.ttk.pte.hu }
\date{}
\maketitle

\centerline{Bull. Math. Soc. Sci. Math. Roumanie 55(103) (2012),
423--428}

\begin{abstract} We prove by using simple number-theoretic
arguments formulae concerning the number of elements of a fixed
order and the number of cyclic subgroups of a direct product of
several finite cyclic groups. We point out that certain
multiplicative properties of related counting functions for finite
Abelian groups are immediate consequences of these formulae.
\end{abstract}

{\sl 2010 Mathematics Subject Classification}: Primary: 20K01,
20K27, Secondary: 11A25.

{\sl Key Words and Phrases}: finite Abelian group, cyclic subgroup,
order of elements, multiplicative arithmetic function.

\section{Introduction}

Consider an arbitrary finite Abelian group $G$ of order $\# G>1$.
Let $c_{\delta}(G)$ and $c(G)$ denote the number of cyclic subgroups
of order $\delta$ ($\delta \mid \# G$) and the number of all cyclic
subgroups of $G$, respectively.

Let $\# G= p_1^{a_1}\cdots p_r^{a_r}$ be the prime power
factorization of $\# G$ and let $G\simeq G_1\times \cdots \times G_r$ be
the primary decomposition of the group $G$, where $\# G_i =
p_i^{a_i}$ ($1\le i\le r$). Furthermore, let $L(G)$ and $L(G_i)$
denote the subgroup lattices of $G$ and $G_i$ ($1\le i\le r$),
respectively. It is known the lattice isomorphism $L(G)\simeq L(G_1)
\times \cdots \times L(G_r)$. See R.~Schmidt \cite{Sch1994},
M.~Suzuki \cite{Suz1951}. It follows that
\begin{equation} \label{c_d_multipl}
c_{\delta}(G)=c_{\delta_1}(G_1)\cdots c_{\delta_r}(G_r),
\end{equation}
where $\delta=\delta_1\cdots \delta_r$ and $\delta_i=p_i^{b_i}$ with
some exponents $b_i\le a_i$ ($1\le i\le r$). Also,
\begin{equation} \label{c_multipl}
c(G)=c(G_1)\cdots c(G_r).
\end{equation}

Consider now a $p$-group $G_{(p)}$ of type
$(\lambda_1,\ldots,\lambda_k)$ with $1\le \lambda_1\le \ldots \le
\lambda_k$. It is known that the number of cyclic subgroups of order
$p^{\nu}$ of $G_{(p)}$ is
\begin{equation} \label{c_p_group}
c_{p^{\nu}}(G_{(p)})= \frac{p^{k-j_{\nu}}-1}{p-1}
p^{\lambda_0+\lambda_1+\ldots+\lambda_{j_{\nu}}+ (k-j_{\nu}-1)(\nu
-1)},
\end{equation}
where $\lambda_{j_{\nu}} <\nu \le \lambda_{j_{\nu}+1}$ and
$\lambda_0=0$. See G.~A.~Miller \cite{Mil1939}, Y.~Yeh
\cite{Yeh1948}.

For every finite Abelian group $G$ the value of $c_{\delta}(G)$ is
therefore completely determined by formulae \eqref{c_d_multipl} and
\eqref{c_p_group}. Obviously, $c(G)=\sum_{\delta \mid \# G}
c_{\delta}(G)$, leading to a formula for $c(G)$.

Recently, formula \eqref{c_p_group} was recovered by
M.~T\u{a}rn\u{a}uceanu \cite[Th.\ 4.3]{Tar2010}, in a slightly
different form, by a method based on properties of certain matrices
attached to the invariant factor decomposition of an Abelian group
$G$, i.e., $G\simeq C_{i_1}\times \cdots \times C_{i_s}$, where
$i_1\cdots i_s=\# G$, $i_1\mid i_2\mid \ldots \mid i_s$, $C_k$
denoting the cyclic group of order $k$.

For arbitrary positive integers $r,n_1,\ldots, n_r$ consider in what
follows the direct product $C_{n_1}\times \cdots \times C_{n_r}$ and
denote by $c(n_1,\ldots,n_r)$ the number of its cyclic subgroups.
Let $\phi$ denote Euler's function.

\begin{theorem} \label{Th_1} For any $n_1,\ldots, n_r\ge 1$, $c(n_1,\ldots,n_r)$
is given by the formula
\begin{equation} \label{total_number_cyclic_subgroups}
c(n_1,\ldots,n_r)= \sum_{d_1\mid n_1,\ldots,d_r\mid n_r}
\frac{\phi(d_1) \cdots \phi(d_r)}{\phi(\lcm(d_1,\ldots,d_r))}.
\end{equation}
\end{theorem}

Formula \eqref{total_number_cyclic_subgroups} was derived in
\cite{Tot2011} using the orbit counting lemma (Burnside's lemma).
Consequently, the number of all cyclic subgroups of a finite Abelian
group can be given by formula \eqref{total_number_cyclic_subgroups}
in a more compact form. Namely, consider the invariant factor
decomposition of $G$ given above and apply
\eqref{total_number_cyclic_subgroups} by selecting $r=s$ and
$n_1=i_1,\ldots,n_s=i_s$. Note that
\eqref{total_number_cyclic_subgroups} can be used also starting with
the primary decomposition of $G$.

In the case $r=2$, by using the identity $\phi(d_1)\phi(d_2)=
\phi(\gcd(d_1,d_2)) \phi(\lcm(d_1,d_2))$ ($d_1,d_2\ge 1$) we deduce
from \eqref{total_number_cyclic_subgroups} the next result:

\begin{cor} For every $n_1,n_2\ge 1$ the number of cyclic
subgroups of $C_{n_1}\times C_{n_2}$ is
\begin{equation} \label{total_number_cyclic_subgr_2_1}
c(n_1,n_2)=  \sum_{d_1 \mid n_1, d_2\mid n_2} \phi(\gcd(d_1,d_2)).
\end{equation}
\end{cor}

It is the purpose of the present paper to give another direct proof
of Theorem \ref{Th_1}, using simple number-theoretic arguments. In
fact, we prove the following formulae concerning the number
$o_{\delta}(n_1,\ldots,n_r)$ of elements of order $\delta$ in
$C_{n_1}\times \cdots \times C_{n_r}$. Let
$n:=\lcm(n_1,\ldots,n_r)$. Obviously, the order of every element of
the direct product is a divisor of $n$. Let $\mu$ be the M\"obius
function.

\begin{theorem} \label{Th_2} For every $n_1,\ldots, n_r\ge 1$ and every $\delta \mid n$,
\begin{equation} \label{number_elements_order_d_first}
o_{\delta}(n_1,\ldots,n_r)= \sum_{e\mid \delta}  \gcd(e,n_1) \cdots
\gcd(e,n_r) \mu(\delta/e)
\end{equation}
\begin{equation} \label{number_elements_order_d_second}
= \sum_{\substack{d_1\mid n_1,\ldots, d_r\mid n_r
\\ \lcm(d_1,\ldots,d_r)=\delta}} \phi(d_1) \cdots \phi(d_r).
\end{equation}
\end{theorem}

Let $c_{\delta}(n_1,\ldots,n_r)$ denote the number of cyclic
subgroups of order $\delta$ ($\delta \mid n$) of the group
$C_{n_1}\times \cdots \times C_{n_r}$. Since a cyclic subgroup of
order $\delta$ has $\phi(\delta)$ generators,
\begin{equation} \label{number_cyclic_subgroups_order_d}
c_{\delta}(n_1,\ldots,n_r)=\frac{o_{\delta}(n_1,\ldots,n_r)}{\phi(\delta)}.
\end{equation}

Now \eqref{total_number_cyclic_subgroups} follows immediately from
\eqref{number_elements_order_d_second} and
\eqref{number_cyclic_subgroups_order_d} by
$c(n_1,\ldots,n_r)=\sum_{\delta\mid n} c_{\delta}(n_1,\ldots,n_r)$.

\begin{cor} For every prime $p$ and every $a_1,\ldots,a_r\ge 1$, $1\le
\nu \le \max(a_1,\ldots,a_r)$,
\begin{equation} \label{number_p}
c_{p^{\nu}}(p^{a_1},\ldots,p^{a_r})= \frac1{p^{\nu-1}(p-1)} \left(
p^{\min(\nu,a_1)+\ldots+\min(\nu,a_r)}-
p^{\min(\nu-1,a_1)+\ldots+\min(\nu-1,a_r)}\right).
\end{equation}
\end{cor}

Here \eqref{number_p} follows at once by
\eqref{number_elements_order_d_first}. In case of a $p$-group of
type $(\lambda_1,\ldots,\lambda_k)$ and choosing $r=k$,
$a_1=\lambda_1,\ldots, a_k=\lambda_k$, \eqref{number_p} reduces to
\eqref{c_p_group}.

The proof of Theorem \ref{Th_2} is contained in Section
\ref{section_2}. Additional remarks are given in Section
\ref{section_3}. We point out, among others, that the multiplicative
properties \eqref{c_d_multipl} and \eqref{c_multipl} are direct
consequences of the formulae included in Theorems \ref{Th_1} and
\ref{Th_2}.


\section{Proof of Theorem \ref{Th_2}} \label{section_2}

Let $o(x)$ denote the order of a group element $x$. Let
$G=C_{n_1}\times \cdots \times C_{n_r}$ in short. Let $x_i\in
C_{n_i}$ such that $o(x_i)=n_i$ ($1\le i\le r$). Then $G= \{
x=(x_1^{i_1},\ldots,x_r^{i_r}): 1\le i_1\le n_1, ..., 1\le i_r\le
n_r\}$.

Using elementary properties of the order of elements in a group and
of the gcd and lcm of integers, respectively we deduce that for
every $x=(x_1^{i_1},\ldots,x_r^{i_r})\in G$,
\begin{equation*}
o(x)= \lcm(o(x_1^{i_1}),\ldots,o(x_r^{i_r}))= \lcm
\left(\frac{o(x_1)}{\gcd(o(x_1),i_1)},\ldots,
\frac{o(x_r)}{\gcd(o(x_r),i_r)} \right)
\end{equation*}
\begin{equation*}
= \lcm \left(\frac{n_1}{\gcd(n_1,i_1)},\ldots,
\frac{n_r}{\gcd(n_r,i_r)} \right) = \lcm \left(\frac{n}{\gcd(n, i_1
n/n_1)},\ldots, \frac{n}{\gcd(n,i_r n/n_r)} \right)
\end{equation*}
\begin{equation*}
= \frac{n}{\gcd (i_1 n/n_1,\ldots, i_r n/n_r,n)}.
\end{equation*}

Assume that $o(x)=\delta$, where $\delta \mid n$ is fixed. Then
$\gcd (i_1 n/n_1,\ldots, i_r n/n_r,n)=n/\delta$.

Write $i_1 n/n_1=j_1 n/\delta, \ldots, i_r n/n_r=j_r n/\delta$. Then
$\gcd(j_1,\ldots,j_r,\delta)=1$ and $j_1=\delta i_1/n_1$, ...,
$j_r=\delta i_r/n_r$ are integers, that is $\delta i_1\equiv 0$ (mod
$n_1$), ..., $\delta i_r\equiv 0$ (mod $n_r$). We obtain, as solutions
of these linear congruences, that
$i_1=k_1n_1/\gcd(\delta,n_1)$ with $1\le k_1\le \gcd(\delta,n_1)$,
..., $i_r=k_rn_r/\gcd(\delta,n_r)$ with $1\le k_r\le
\gcd(\delta,n_r)$.

Therefore, the number of elements of order $\delta$ in $G$ is the
number of ordered $r$-tuples $(k_1,\ldots,k_r)$ satisfying the
conditions

(i)  $1\le k_1\le \gcd(\delta,n_1),\ldots, 1\le k_r\le
\gcd(\delta,n_r)$ and

(ii) $\gcd \left(k_1\frac{\delta}{\gcd(\delta,n_1)},\ldots,
k_r\frac{\delta}{\gcd(\delta,n_r)}, \delta \right)=1$.

Using the familiar formula $\sum_{d\mid n} \mu(d)=1$ for $n=1$ and
$0$ otherwise, this can be written as
\begin{equation*}
o_{\delta}(n_1,\ldots,n_r) = \sum_{e\mid \delta} \mu(e) \sum_{\substack{k_1\le \gcd(\delta,n_1)\\
e\mid k_1\frac{\delta}{\gcd(\delta,n_1)}}} \cdots
\sum_{\substack{k_r\le \gcd(\delta,n_r)\\ e\mid
k_r\frac{\delta}{\gcd(\delta,n_r)}}} 1.
\end{equation*}

Here the linear congruence $k_1\frac{\delta}{\gcd(\delta,n_1)}
\equiv 0$ (mod $e$) has $\gcd(e,\frac{\delta}{\gcd(\delta,n_1)})$
solutions in $k_1$ (mod $e$), and (mod $\gcd(\delta,n_1)$) it has
exactly
\begin{equation*}
\frac{\gcd(\delta,n_1)}{e}
\gcd\left(e,\frac{\delta}{\gcd(\delta,n_1)}\right)=
\gcd\left(\gcd(\delta,n_1),\frac{\delta}{e}\right) =
\gcd\left(\frac{\delta}{e},n_1\right)
\end{equation*}
solutions (similar for the indexes $2,...,r$). Hence,
\begin{equation*}
o_{\delta}(n_1,\ldots,n_r) = \sum_{e\mid \delta} \mu(e)
\gcd(\delta/e,n_1) \cdots \gcd(\delta/e,n_r)
\end{equation*}
\begin{equation*}
= \sum_{e\mid \delta} \mu(\delta/e) \gcd(e,n_1) \cdots \gcd(e,n_r),
\end{equation*}
which is formula \eqref{number_elements_order_d_first}.

Now using the identity $\sum_{d\mid n} \phi(d)=n$,
\begin{equation*}
o_{\delta}(n_1,\ldots,n_r) = \sum_{e\mid \delta} \mu(\delta/e)
\sum_{d_1\mid \gcd(e,n_1)} \phi(d_1) \cdots \sum_{d_r\mid
\gcd(e,n_r)} \phi(d_r)
\end{equation*}
\begin{equation*}
= \sum_{d_1\mid n_1,\ldots, d_r\mid n_r} \phi(d_1)\cdots \phi(d_r)
\sum_{ab \lcm(d_1,\ldots,d_r)=\delta} \mu(a),
\end{equation*}
where the inner sum is $0$ unless $\lcm(d_1,\ldots,d_r)=\delta$, and
in this case it is $1$. This gives
\eqref{number_elements_order_d_second}. The proof of Theorem
\ref{Th_2} is complete.


\section{Remarks}  \label{section_3}

The function $c(n_1,\ldots,n_r)$ is representing a multiplicative
function of $r$ variables, i.e.,
\begin{equation} \label{c_multiplic}
c(n_1m_1,\ldots,n_rm_r)= c(n_1,\ldots,n_r) c(m_1,\ldots,m_r)
\end{equation}
holds for any $n_1,\ldots,n_r,m_1,\ldots,m_r\ge 1$ such that
$\gcd(n_1\cdots n_r,m_1\cdots m_r)=1$. Therefore,
\begin{equation*}
c(n_1,\ldots,n_r)= \prod_p c(p^{e_p(n_1)}, \ldots,p^{e_p(n_r)}),
\end{equation*}
where $n_i=\prod_p p^{e_p(n_i)}$ is the prime power factorization of
$n_i$, the product being over the primes $p$ and all but a finite
number of the exponents $e_p(n_i)$ are zero ($1\le i\le r$).

The property \eqref{c_multiplic} follows by the next simple
number-theoretic argument: According to
\eqref{total_number_cyclic_subgroups}, $c(n_1,\ldots,n_r)$ is the
$r$ variables convolution of the function $(n_1,\ldots,n_r) \mapsto
\frac{\phi(n_1) \cdots \phi(n_r)}{\phi(\lcm(n_1,\ldots,n_r))}$ with
the constant $1$ function, both multiplicative. Since convolution
preserves the multiplicativity we deduce that $c(n_1,\ldots,n_r)$ is
multiplicative, too. See \cite[Section 2]{Tot2011} for details.

Also, let $n_1,\ldots,n_r,m_1,\ldots,m_r\ge 1$ such that
$\gcd(n_1\cdots n_r,m_1\cdots m_r)=1$. Let
$n:=\lcm(n_1,\ldots,n_r)$, $m:=\lcm(m_1,\ldots,m_r)$. For $\delta
\mid nm$ write $\delta=ab$ with $a\mid n$, $b\mid m$ (which are
unique). Then,
\begin{equation} \label{c_delta_multipl}
c_{\delta}(n_1m_1,\ldots,n_rm_r)= c_a(n_1,\ldots,n_r)
c_b(m_1,\ldots,m_r),
\end{equation}
which follows from \eqref{number_elements_order_d_first} by a short
computation.

Now, \eqref{c_multiplic} and \eqref{c_delta_multipl} show the
validity of \eqref{c_multipl} and \eqref{c_d_multipl} quoted in the
Introduction.

Let $A(n_1,\ldots,n_r)$ denote the arithmetic mean of the orders of
elements in $C_{n_1}\times \cdots \times C_{n_r}$. From
\eqref{number_elements_order_d_second} we deduce

\begin{cor} For every $n_1,\ldots,n_r\ge 1$,
\begin{equation*}
A(n_1,\ldots,n_r)= \frac1{n_1\cdots n_r} \sum_{d_1\mid n_1,\ldots,
d_r\mid n_r} \phi(d_1)\cdots \phi(d_r) \lcm(d_1,\ldots,d_r).
\end{equation*}
\end{cor}

The function $(n_1,\ldots,n_r) \mapsto A(n_1,\ldots,n_r)$ is also
multiplicative. See \cite{zGat2004,LucShp2003} for further
properties of the average of orders of elements in finite Abelian
groups.

Consider the special case $n_1=\ldots =n_r=n$. According to
\eqref{number_elements_order_d_first},
\begin{equation} \label{Jordan}
o_{\delta}(n,\ldots,n) = \sum_{e\mid \delta} e^r \mu(\delta/r) =
\phi_r(\delta)
\end{equation}
is the Jordan function of order $r$. Note that
$\phi_1(\delta)=\phi(\delta)$ is Euler's function.

This means that the number of elements of order $\delta$ ($\delta
\mid n$) in the direct product $C_n\times \cdots \times C_n$, with
$r$ factors, is $\phi_r(\delta)= \delta^r\prod_{p\mid \delta}
(1-1/p^r)$. This is a less known group-theoretic interpretation of
the Jordan function. See G.~A.~Miller \cite{Mil1904,Mil1905},
J.~Schulte \cite[Th.\ 7]{Sch1999}.

Hence the number of cyclic subgroups of order $\delta$ ($\delta
\mid n$) of the group $C_n\times \cdots \times C_n$, with
$r$ factors, is $\phi_r(\delta)/\phi(\delta)$.

Note that the function $k\mapsto \phi_r(k)/(k^{r-2}\phi(k))=
k\prod_{p\mid k} (1+1/p+1/p^2+\ldots+1/p^{r-1})$ ($r\ge 2$) was
considered quite recently in \cite{SolPla2011}, in connection with
Robin's inequality and the Riemann hypothesis.

As a byproduct we obtain from \eqref{number_elements_order_d_second}
and \eqref{Jordan} that for every positive integer $\delta$,
\begin{equation} \label{Sterneck}
\sum_{\lcm(d_1,\ldots,d_r)=\delta} \phi(d_1)\cdots \phi(d_r)
=\phi_r(\delta),
\end{equation}
where the sum is over all ordered $r$-tuples $(d_1,\ldots,d_r)$ such
that $\lcm(d_1,\ldots,d_r)=\delta$. The function in the left hand
side of \eqref{Sterneck} is the von Sterneck function, cf. \cite[p.\
14]{McC1986}.

A final remark: The number $s(n_1,n_2)$ of all subgroups of the
group $C_{n_1}\times C_{n_2}$ is given by
\begin{equation} \label{s_n}
s(n_1,n_2)= \sum_{d_1\mid n_1, d_2\mid n_2} \gcd(d_1,d_2) \quad
(n_1,n_2\ge 1).
\end{equation}

Compare \eqref{s_n} to \eqref{total_number_cyclic_subgr_2_1}. A
proof of \eqref{s_n} will be presented elsewhere.


\end{document}